\documentclass[10pt]{amsart}

\usepackage{amscd}
\usepackage{eucal}
\usepackage{enumerate}
\usepackage{amsthm}
\usepackage{amssymb} 
\usepackage{amsmath}
\usepackage[dvips]{graphicx}
\usepackage{psfrag}
\usepackage[all]{xy}
\usepackage{epsfig}

\newtheorem{Theorem}{Theorem}[section]
\newtheorem*{theorem*}{Theorem 1}

\newtheorem*{proposition*}{Proposition}
\newtheorem{proposition}{Proposition}[section]

\newtheorem*{corollary*}{Corollary}
\theoremstyle{definition}
\newtheorem*{definition}{Definition}
\newtheorem*{comments*}{Comments}
\newtheorem{example}{Example}
\newtheorem*{example*}{Example}

\newtheorem{remark}{Remark}
\newtheorem*{remarks*}{Remarks}

\numberwithin{equation}{section}

\gdef\myletter{}

\let\savetheequation\theequation
\def\theequation{\savetheequation\myletter}

\newcommand{\bbC}{{\mathbb C}}

\newcommand{\bbZ}{{\mathbb Z}}
\newcommand{\bbQ}{{\mathbb Q}}

\newcommand{\m}{{\mathrm{min}}}
\newcommand{\M}{{\mathrm{max}}}
\newcommand{\vol}{{\mathrm{vol}}}

\newcommand{\om}{{\omega}}

\newcommand{\Ga}{{\Gamma}}

\newcommand{\la}{{\lambda}}

\newcommand{\bN}{{\bf N}}
\newcommand{\hN}{{\hat{N}}}

\newcommand{\T}{{\mathbb T}}

\newcommand{\Z}{{\mathbb Z}}
\newcommand{\C}{{\mathbb C}}

\newcommand{\R}{{\mathbb R}}

\newcommand{\scal}{{\tau}}

\def    \to     {{\longrightarrow}}

\begin{document}
\title{Hearing the weights of weighted projective planes}
\author{Miguel Abreu} 
\thanks{M. Abreu and L. Godinho were partially supported by FCT through grant POCTI/MAT/57888/2004, and L. Godinho was also partially supported by the Funda\c{c}\~{a}o Calouste Gulbenkian; P. Freitas was partially supported by FCT through grant POCTI/0208/2003.}
\address{Departamento de Matem\'{a}tica, Instituto Superior T\'{e}cnico, Av. Rovisco Pais, 1049-001 Lisboa, Portugal}
\email{mabreu@math.ist.utl.pt}
\author{Emily B. Dryden}
\address{Centro de An\'{a}lise Matem\'{a}tica, Geometria e Sistemas Din\^{a}micos, Instituto Superior T\'{e}cnico, Av. Rovisco Pais, 1049-001 Lisboa, Portugal}
\email{dryden@math.ist.utl.pt}
\author{Pedro Freitas}
\address{Departamento de Matem\'{a}tica, Faculdade de Motricidade Humana (TU Lisbon) {\rm and}
Mathematical Physics Group of the University of Lisbon, Complexo Interdisciplinar, Av. Prof. Gama Pinto 2, P-1649-003 Lisboa, Portugal}
\email{freitas@cii.fc.ul.pt}
\author{Leonor Godinho}
\address{Departamento de Matem\'{a}tica, Instituto Superior T\'{e}cnico, Av. Rovisco Pais, 1049-001 Lisboa, Portugal}
\email{lgodin@math.ist.utl.pt}

\begin{abstract}
Which properties of an orbifold can we ``hear,'' i.e., which topological and geometric properties of an orbifold are determined by its Laplace spectrum?
We consider this question for a class of four-dimensional K\"{a}hler orbifolds: weighted projective planes $M:=\C P^2(N_1,N_2,N_3)$ with three isolated singularities.  We show that the spectra of the Laplacian acting on $0$- and $1$-forms on $M$ determine the weights $N_1$, $N_2$, and $N_3$.  The proof involves analysis of the heat invariants using several techniques, including localization in equivariant cohomology.  We show that we can replace knowledge of the spectrum on $1$-forms by knowledge of the Euler characteristic and obtain the same result.  Finally, after determining the values of $N_1$, $N_2$, and $N_3$, we can hear whether $M$ is endowed with an extremal K\"{a}hler metric.
\end{abstract}

\subjclass[2000]{Primary 58J50 Secondary 53D20, 55N91} 
\keywords{Laplace spectrum, heat invariants, weighted projective planes}

\maketitle
\section{Introduction}
\label{sec:1}

An orbifold consists of a Hausdorff topological space together with an atlas of
coordinate charts satisfying certain equivariance conditions (cf. \S
\ref{sec:2.1}).  We will be interested in compact K\"{a}hler orbifolds; 
in particular, we focus on weighted projective planes, asking what the Laplace spectrum 
of such a space can tell us about its properties.  

For manifolds, the asymptotic expansion of the heat kernel can be used to
connect the geometry of the manifold to its spectrum.  The so-called \emph{heat
invariants}  appearing in the asymptotic expansion tell us the dimension, 
volume, and various quantities involving the curvature of the manifold.  
We are interested in the K\"ahler setting; using the heat invariants, results have been 
obtained on spectral determination of complex projective spaces \cite{BGM}, Einstein 
manifolds \cite{Sakai, Ta2}, and general K\"ahler manifolds  \cite{Ta1, Ta2}.  These results 
are all for smooth manifolds, while we will be interested in the effects of the presence 
of singularities. 

Orbifolds began appearing sporadically in the spectral theory literature in the
early 1990's, and have been receiving more dedicated attention in the last five
years.  It is known that the volume and dimension of an orbifold are spectrally determined \cite{Far01}.  
The nature of the singularities allowed to appear in an
isospectral set has been studied; in general, 
there can be at most finitely many isotropy types (up to isomorphism)
in a set of isospectral Riemannian orbifolds that share a uniform lower bound on
Ricci curvature \cite{Sta05}.  On the other hand, N. Shams, E. Stanhope and D. Webb
\cite {SSW} constructed arbitrarily large (but always finite) isospectral sets
which satisfy this curvature condition, where each element in a given set has
points of distinct isotropy.  One interesting question is whether the spectrum
determines the orders of the singular points for large classes of orbifolds; we
show that this is the case for weighted projective planes.  In particular, we prove 
(see \S \ref{sec:mainthm})

\begin{theorem*}\label{thm:1}
Let $M:=\C P^2(N_1,N_2,N_3)$ be a four-dimensional weighted projective space with isolated singularities, equipped with a K\"{a}hler metric.  Then the spectra of the Laplacian acting on $0$- and $1$-forms on $M$ determine the weights $N_1, N_2, N_3$.
\end{theorem*}  

\noindent Note that we need to consider the spectrum of the Laplacian acting on both $0$- and $1$-forms.  We conjecture that the spectrum of the Laplacian acting on $0$-forms determines the weights, and have verified this with extensive computer testing.  However, we cannot prove Theorem \ref{thm:1} using only this information at this point (cf. Section \ref{sec:01}) 

The proof  of Theorem \ref{thm:1} involves introducing new techniques into the analysis of the heat invariants.  In particular, we use localization in equivariant cohomology to obtain expressions for topological invariants of $\C P^2(N_1,N_2,N_3)$ in terms of $N_1$, $N_2$, and $N_3$.  We also use a decomposition of the Riemannian curvature tensor inside the space of K\"ahlerian curvature tensors to exploit the complex structure on our weighted projective spaces.  We show that we can replace knowledge of the spectrum on $1$-forms by knowledge of the Euler characteristic and obtain the same result.  Finally, after determining the values of $N_1$, $N_2$, and $N_3$, we can hear whether the metric on $\C P^2(N_1,N_2,N_3)$ is an extremal K\"{a}hler metric in the sense of Calabi \cite{Ca}.

The same methods can be used with higher-dimensional weighted projective spaces to
obtain conditions on the weights (cf. Remark~\ref{rmk:3} in \S \ref{sec:mainthm}). However, these
conditions do not uniquely determine the weights. To obtain additional information we
would need to impose restrictions on the metric and to use higher-order terms of the
asymptotic expansion of the heat trace.

The paper is organized as follows.  We begin with the necessary background on orbifolds and orbi-bundles, then introduce the localization formula for orbifolds and give the setup necessary to apply it in our case.  The decomposition of the Riemannian curvature in terms of K\"ahlerian curvature tensors is presented.  Localization is then used to compute various integrals whose values we will need in the proof of Theorem \ref{thm:1}.  In \S \ref{sec:A}, we use a certain polytope associated to a given symplectic toric orbifold to give a description of extremal metrics on weighted projective planes, and to calculate the integral of the square of the scalar curvature on a weighted projective plane.  In \S \ref{sec:heat}, we recall the asymptotic expansion of the heat trace of an orbifold as given in \cite{DGGWZ}, and we calculate the first few heat invariants for $\C P^2(N_1,N_2,N_3)$.  Finally, we bring all of these tools together to prove Theorem \ref{thm:1} and related results in \S \ref{sec:mainthm}.   



\section{Preliminaries}
\label{sec:2}
\subsection{Orbifolds and Orbi-bundles}
\label{sec:2.1}
An \emph{orbifold} $M$ is a singular manifold whose singularities are locally isomorphic to quotient singularities of the form $\R^n/\Gamma$, where $\Gamma$ is a finite subgroup of $GL(n,\R)$. An \emph{orbifold chart} on $M$ is a triple $(U,\Gamma, V)$ consisting of an open subset $U$ of $M$, a finite group $\Gamma \subset GL(n,\R)$, an open subset $V$ of $\R^n$ and a homeomorphism $U=V/\Gamma$. An orbifold structure on the paracompact Hausdorff space $|M|$ is then a collection of orbifold charts $\{(U_i, \Gamma_i, V_i)\}$ covering $|M|$, subject to appropriate compatibility conditions. In particular, these conditions ensure that on each connected component of $M$, the generic stabilizers of the $\Gamma_i$-actions are isomorphic.  For each singularity $p\in M$, there is a finite subgroup $\Gamma_p$ of $GL(n,\R)$, unique up to conjugation, such that open neighborhoods of $p$ in $M$ are homeomorphic to neighborhoods of the origin in $\R^n/\Gamma_p$.  We call $\Gamma_p$ the \emph{orbifold structure group} of $p$. 
\begin{example}
\label{example:0}

Let $\bN = (N_1, \ldots, N_{m+1})$ be a vector of positive integers
which are pairwise relatively prime. The weighted projective space
\[
\C P^{m}(\bN) := \C P^{m}(N_1, \ldots, N_{m+1}) :=
(\C^{m+1})^\ast / \sim\,,
\]
where
\[
((z_1,\ldots,z_{m+1}) \sim (\lambda^{N_1} z_1, 
\ldots, \lambda^{N_{m+1}} z_{m+1}),\, \lambda\in \C^*)\,,
\]
is a compact orbifold.  It has $m+1$ isolated singularities at the points 
$[1:0: \cdots:0],\,\ldots,\, [0:\cdots:0:1]$, with
orbifold structure groups $\Z_{N_1},\,\ldots,\,\Z_{N_{m+1}}$. Note
that $\C P^m({\bf 1})$ is the usual smooth projective space $\C P^m$.
\end{example}

The definitions of vector fields, differential forms, metrics, group actions, etc. naturally extend to orbifolds. For instance, a \emph{symplectic orbifold} is an orbifold equipped with a closed non-degenerate $2$-form $\omega$. 

Given any orbifold $M$, \emph{fiber orbi-bundles} (or \emph{orbifold bundles}) $\pi:E\to M$ are defined by $\Gamma$-equivariant fiber bundles $Z\to E_V\to V$ in orbifold charts $(U,\Gamma, V)$, together with suitable compatibility conditions. Each fiber $\pi^{-1}(p)$ is not, in general, diffeomorphic to $Z$, but to some quotient of $Z$ by the action of the orbifold structure group $\Gamma_p$.  

Line orbi-bundles $\pi: L \to \Sigma$ over a $2$-dimensional orbifold $\Sigma$ (also called an \emph{orbi-surface}) of genus $g$ with $k$ cone singularities with orders $N_1, \ldots, N_k$, correspond to Seifert fibrations ($3$-dimensional manifolds together with an $S^1$ action with finite stabilizers \cite{Se, S, FS}) when we take the corresponding circle bundles. Hence, to each line orbi-bundle we can associate a collection of integers $(b, m_1, \cdots, m_k)$, called the \emph{Seifert invariant} of $L$ over $\Sigma$,  as well as its \emph{Chern number} or  \emph{degree} defined by the formula $\mathrm{deg}(L):=b+\sum_{i=1}^k \frac{m_i}{N_i}$.  Let us briefly recall how to obtain this invariant (see for example \cite{MOY} and \cite{FS} for a detailed construction).
If $x_1, \ldots, x_k$ are the orbifold singularities of $\Sigma$ and $N_1,\ldots, N_k$ are the orders of the corresponding orbifold structure groups, then a neighborhood of $x_i$ in $L$ is of the form $(D^2 \times \bbC)/\bbZ_{N_i}$, where $D^2$ is a $2$-disk and where $\bbZ_{N_i}$ acts on $D^2 \times \bbC$ by
$$
\xi_{N_i}\cdot (z,w)=(\xi_{N_i} z, \xi_{N_i}^{m_i} w), 
$$ 
for some integer $0 \leq m_i <N_i$, and $\xi_{N_i}$ a primitive $N_i$-th root of unity.
Given $\Sigma$ and its orbifold singularities $x_1,\ldots, x_k$, we define special line orbi-bundles $H_{x_i}$ as follows.  Let $\Sigma_{x_i}$ be the orbi-surface obtained from $\Sigma$ by deleting a small open neighborhood $U_i$ around $x_i$, where $U_i$ is isomorphic to $D^2/\bbZ_{N_i}$.  Take the trivial bundle over $\Sigma_{x_i}$, and over $U_i$ take the line orbi-bundle $(D^2\times \bbC)/\bbZ_{N_i}$, for a $\bbZ_{N_i}$-action given by
$$
\xi_{N_i}\cdot (z,w)=(\xi_{N_i} z, \xi_{N_i} w)
$$ 
with gluing map $\alpha : \partial \Sigma_{x_i} \times \bbC \to  (\partial D^2\times \bbC)/\bbZ_{N_i}$, 
$$
\alpha(e^{i\theta},w)=(e^{-i\theta},e^{-i\theta} w).
$$
That is, 
$$
H_{x_i} := \left(\Sigma^+_{x_i} \times \bbC \right) \bigcup_\alpha \, (D^{2^-}\times \bbC)/\bbZ_{N_i},
$$
where $\Sigma^+_{x_i}$ is positively oriented and $D^{2^-}$ is negatively oriented.
The bundle $L \otimes H_{x_1}^{-m_1} \otimes \cdots \otimes H_{x_k}^{-m_k}$, called the \emph{de-singularization} of $L$,  is a trivial line orbi-bundle over each neighborhood of the $x_i$'s and is naturally isomorphic to a smooth line bundle $|L|$ over the (smooth) surface $|\Sigma|$.  Taking the first Chern number $b$ of $|L|$, the collection of integers $(b,m_1,\ldots, m_k)$  is then the \emph{Seifert invariant} of $L$ over $\Sigma$. 

This invariant of an orbifold line bundle classifies it. Indeed, denoting by $\mathrm{Pic}^t(\Sigma)$ the \emph{Picard group} of topological isomorphism classes of line orbi-bundles over $\Sigma$ (with group operation the tensor product),  the map 
\begin{align*}
\mathrm{Pic}^t(\Sigma) & \,\,\to \,\, \bbQ \oplus \bigoplus_{i=1}^k \bbZ_{N_i} \\
 L & \mapsto (\mathrm{deg}(L), m_1,\ldots, m_k),
\end{align*} 
is an injection with image the set of tuples $(c,  m_1,\ldots, m_k)$ with $c = \sum_{i=1}^k m_i/ N_i \pmod{\bbZ}$ (cf. \cite{FS}). In particular, if the $N_i$'s are pairwise relatively prime, then   $\mathrm{Pic}^t(\Sigma) \cong \bbZ$ is generated by a line orbi-bundle $L_0$ with $\mathrm{deg}(L_0)=\frac{1}{N_1\cdots N_k}$.

An important class of examples of line orbi-bundles are those over weighted projective spaces $\Sigma=\bbC P^1(p,q)$. Just as any line bundle over a sphere is isomorphic to some line bundle $\mathcal{O}(r)$ with first Chern number $r\in \Z$, defined by $(S^3\times \C)/\!\!\sim$, where $(z_1,z_2,w) \sim (\lambda z_1, \lambda z_2, \lambda^r w)$ for $\lambda \in S^1$, any  line orbi-bundle with isolated singularities over an orbifold ``sphere'' $\C P^1(p,q)$ with $p$ and $q$ relatively prime  is isomorphic to some $\mathcal{O}_{p,q}(r)=L_0^r$, defined by  $(S^3\times \C)/\!\!\sim$, where now $(z_1,z_2,w) \sim (\lambda^{p} z_1, \lambda^{q} z_2, \lambda^r w)$ for $\lambda \in S^1$. Note that $\mathrm{deg}(\mathcal{O}_{p,q}(r))=\frac{r}{pq}$ and that, if $r>0$, $\mathcal{O}_{p,q}(r)$ is isomorphic to the normal orbi-bundle of $\C P^1(p,q)$ inside the weighted projective space $\C P^2(p, q , r)$.

\subsection{Group actions and equivariant cohomology}
\label{sec:2.2}
Let $(M,\omega)$ be a symplectic orbifold, and $G$ a compact, connected Lie group acting smoothly on $M$. Since $G$ is connected, the components of the fixed point set $M^G$ are suborbifolds of $M$. Moreover, if $G$ is abelian, the normal orbi-bundle $\nu_F$ of each fixed point component $F$ is even dimensional, and admits an invariant Hermitian structure. If $M$ is oriented, any choice of  such a Hermitian structure defines an orientation in $F$. The definition of \emph{equivariant differential forms}, $\mathcal{A}_G(M)$ (that is, polynomial $G$-equivariant mappings $\alpha$ from the Lie algebra $\frak{g}$ of $G$ to the space of differential forms in $M$) extends naturally to orbifolds equipped with a group action, as does the equivariant differential $d_G:\mathcal{A}_G(M)\to\mathcal{A}_G(M)$ defined as
$$
d_G(\alpha)(\xi)=d\alpha(\xi) -2\pi\,i\,\iota(\xi_M)\alpha(\xi),
$$ 
where $\xi_M$ is the fundamental vector field corresponding to $\xi$.

Let us assume now that $G=T$ is abelian, that $M$ is compact, connected and oriented and consider the integration mapping $\int:\mathcal{A}_T(M)\to \mathcal{A}_T(M)$ and  the embeddings $\iota_F:F\to M$ of the connected components of the fixed point set. The Atiyah-Bott and Berline-Vergne localization formula (\cite{AB,BV}), generalized to orbifolds by Meinrenken \cite{M}, states that
\begin{proposition}(Localization formula for orbifolds)
\label{prop:2.1}
Suppose $G=T$ is abelian, and let $\alpha\in \mathcal{A}_T(M)$ be $d_T$-closed. Then
$$
\frac{1}{d_M}\int_M\alpha=\sum_F \frac{1}{d_F}\int_F \frac{\iota_F^*\alpha}{e(\nu_F)}, 
$$
where the sum is over the connected components of the fixed point set, where, for a connected orbifold $X$, $d_X$ is the order of the orbifold structure group of a generic point of $X$, and where $e(\nu_F)$ is the equivariant Euler class of  $\nu_F$, the normal orbi-bundle to $F$.
\end{proposition}

The right hand side of this equation is very simple if $M$ admits an invariant almost complex structure. In this case, the computation of the equivariant Chern classes of $\nu_F$ is given by  the corresponding \emph{equivariant Chern series} $c^{T}(\nu_F):=\sum_i t^i\, c_i^{T}(\nu_F)$. Moreover (using the splitting principle if necessary) we can assume, without loss of generality, that $\nu_F$ splits into a direct sum of invariant line orbi-bundles $L_i$ with first Chern classes $c_1(L_i)$  where $T$ acts with rational \emph{orbi-weights}\footnote{Given an orbifold chart $(U,\Gamma, V)$ around a point in $F$ it is not always true that the $G$-action on $U$ lifts to $V$ but some finite covering $\hat{G} \to G$ does. The weights for this action of $\hat{G}$ on $\nu_F$ are called the \emph{orbi-weights} of $G$.} $\lambda_i$, and we have:
\begin{equation} 
\label{eq:2.0}
c_t^{T}(\nu_F)=\prod_i\left(1+t\,(c_1(L_i)+\lambda_i) \right).
\end{equation}

As an example, the \emph{equivariant Euler class} $e(\nu_F)$  is
\begin{equation} 
\label{eq:2.1}
e(\nu_F)=\prod_i\left(c_1(L_i)+\lambda_i \right)
\end{equation}
and the \emph{first equivariant Chern class} $c_1^{T}(\nu_F)$ is 
\begin{equation} 
\label{eq:2.2}
c_1^{T}(\nu_F)=\sum_i \left(c_1(L_i)+\lambda_i \right).
\end{equation}

\subsection{Circle actions}
\label{sec:2.3}
Let $M$ be a $4$-dimensional symplectic orbifold with isolated cone singularities equipped with a Hamiltonian $S^1$-action, 
and let $F$ be a fixed point. A neighborhood of $F$ is modeled by some quotient $\C^2/\Z_N$ for a $\Z_N$ action given by
$$
\xi_N\cdot(z,w)=(\xi_N z,\xi_N^{m}w),
$$     
with $1 \leq m < N$   and $(m,N)=1$  (we are assuming that the orbifold singularities are isolated). The circle action on this neighborhood will be given by
$$
e^{i\textsl{x}}\cdot(z,w) =(e^{i\frac{k_1}{N}\textsl{x}}\,z,e^{i\frac{k_2}{N}\textsl{x}}\,w),
$$
for some integers $k_1$ and $k_2$ (see \cite{Go} for details). If $F$ is an isolated fixed point then the greatest common divisor $(k_1,k_2)$ is equal to $1$ or $N$ (we are assuming the action to be effective), and $k_2=mk_1 \pmod{N}$ (in order to have  a well-defined action). The numbers $\frac{k_1}{N}\textsl{x}$ and $\frac{k_2}{N}\textsl{x}$ are  the orbi-weights of the action at $F$ and, in these coordinates, the moment map is given by $\phi(z,w)=\phi(F) + \frac{1}{2}(\frac{k_1}{N}|z|^2 + \frac{k_2}{N} |w|^2)$. If $F$ is not an isolated fixed point, then the orbi-weight tangent to the fixed surface containing  $F$  is equal to zero while  the normal one  is equal to $\pm \textsl{x}$ (again for the action to be effective). 

\begin{example}\label{example:1} Consider Example~\ref{example:0} with $n=2$, i.e., 
let $M$ be the weighted projective space $\bbC P^2(N_1,N_2,N_3)$ where the positive integers $N_i$ are pairwise relatively prime, now equipped with  the $S^1$-action given by
$$
e^{i\textsl{x}} \cdot [z_0: z_1: z_2] = [z_0: z_1:  e^{i \textsl{x}} z_2].
$$
This action is Hamiltonian with respect to the standard symplectic form $\omega$. Moreover, it fixes the point $F_3:=[0:0:1]$ as well as  the orbi-surface $\Sigma:=\bbC P^1(N_1,N_2)$.  
A coordinate system centered at $F_3$ is given by 
$$
(z,w)\in \C^2 \longmapsto [z:w:1]\in M,
$$ 
with $(z,w) \sim ( \xi_{N_3}^{N_1} \, z , \xi_{N_3}^{N_2}\, w)$, where $\xi_{N_3}$ is  a primitive $N_3$-th root of unity. The action of $S^1$ on this coordinate system is given by
$$
e^{\, i \textsl{x}} \cdot  [z: w: 1]  = [z : w : e^{i \textsl{x}}] = [e^{-i\frac{N_1}{N_3}\textsl{x}} z : e^{-i\frac{N_2}{N_3}\textsl{x}} w: 1],
$$
implying that the weights of the isotropy representation of $S^1$ on $T_{F_3}M$ are $- \frac{N_1}{N_3}\textsl{x}$ and  $- \frac{N_2}{N_3} \textsl{x}$.
On the other hand, the orbi-weights at every point in $\Sigma$ are $(\textsl{x},0)$. Note that the degree of the normal orbi-bundle of $\Sigma$ inside $M$ is $b_\Sigma = \frac{N_3}{N_1\, N_2}$. 
\end{example}

\subsection{Riemann curvature tensor}
\label{sec:2.4}
We now review the well-known decomposition of the Riemann curvature tensor of a Riemannian metric.  A detailed exposition for the manifold case can be found for instance in \cite{B, Gau, A-C-G}.  Since we are interested in weighted projective spaces, which admit a K\"ahler metric, our goal is to connect this standard decomposition of the Riemannian curvature to its K\"ahlerian decomposition.

In what follows, we extract the necessary background from \cite{A-C-G}. The curvature $R$ on a Riemannian orbifold $(M,g)$ is defined, as usual, by 
$$
R_{X,Y} Z = \nabla_{\left[X,Y\right]} Z - \left[\nabla_X,\nabla_Y\right]Z
$$
for all vector fields $X,Y,Z$, where $\nabla$ is the Levi-Civita connection. It is a $2$-form with values in the adjoint orbi-bundle $AM$ (the bundle of skew-symmetric endomorphisms of the tangent orbi-bundle) and satisfies the Bianchi identity\footnote{$R_{X,Y}Z+R_{Y,Z}X+R_{Z,X}Y=0$.}. Using the metric $g$ we can identify $AM$ with $\Omega^2M$ and so $R$ can be viewed as a section of $\Omega^2M \otimes \Omega^2M$. Moreover, it follows from the Bianchi identity  that $R$ belongs to the symmetric part $S^2\Omega^2M$ of $\Omega^2M \otimes \Omega^2M$ as well as to the kernel of the map 
$$
\beta: S^2 \Omega^2 M \to \Omega^4 M
$$  
determined by the wedge product. The kernel of $\beta$, $\mathcal{R}M$, is called the orbi-bundle of abstract curvature tensors.

In addition, we have the \emph{Ricci contraction}, a linear map 
$$ 
c: \mathcal{R}M \to SM,
$$
where $SM$ denotes the orbi-bundle of symmetric bilinear forms on $M$, which sends $R$ to the Ricci form $Ric$.\footnote{$Ric_{X,Y}=\text{tr}(Z \mapsto R_{X,Z} Y)$.} Hence, we obtain an orthogonal decomposition
$$
\mathcal{R}M = c^*(SM) \oplus \mathcal{W}M,
$$
where $\mathcal{W}M$, called the orbi-bundle of \emph{abstract Weyl tensors}, is the kernel of $c$ in $\mathcal{R}M$. Consequently, $R$ can be written as 
$R=c^*(h)+W$, where $W$ is called the \emph{Weyl tensor} and $h$ is such that $c\circ c^*(h)=Ric$. Since, for $n\geq 3$ the map $c$ is surjective, $c^*$ is injective and $h$ is given by
$$
h= \frac{\scal}{2n(n-1)}g + \frac{Ric_0}{n-2},
$$ 
where $\scal$ is the scalar curvature (the trace of $Ric$ with respect to $g$), and where $Ric_0$ denotes the traceless part of $Ric$ (i.e. $Ric=\frac{\scal}{n}g+Ric_0$).
Hence, the curvature $R$, viewed as a symmetric endomorphism of $\Omega^2M$ using $g$, decomposes as
\begin{equation}
\label{eq:decomp}
R=\frac{\scal}{n(n-1)}Id_{\vert_{\Omega^2 M}} + \frac{1}{n-2}\{Ric_0,\cdot\} + W,
\end{equation}
where $\{Ric_0,\cdot\}$ acts on $\alpha \in \Omega^2M$ as the anti-commutator, $\{Ric_0,\alpha \}:= Ric_0\circ \alpha + \alpha \circ Ric_0$.
We will use the standard notation 
\begin{align*}
U:= & \frac{\scal}{n(n-1)}Id_{\vert_{\Omega^2 M}}\\
Z:= & \frac{1}{n-2}\{Ric_0,\cdot\}
\end{align*}
for the first two terms of this orthogonal decomposition (note that in this notation  $|R|^2=|U|^2+|Z|^2+|W|^2$).

Since we are interested in K\"ahler orbifolds, we can consider an orthogonal complex structure $J$ parallel with respect to the Levi-Civita connection and obtain a K\"ahler form $\omega$ as $\omega(X,Y)=g(JX,Y)$. We can also define the \emph{Ricci form} and its primitive part $\rho_0$ (that is, such that $(\rho_0,\omega)=0$) in the same way\footnote{$\rho(X,Y)=Ric(JX,Y)$ and $\rho_0(X,Y)=Ric_0(JX,Y)$}. Note that this isomorphism  $S\mapsto S(J\cdot,\cdot)=:S\circ J$ from the $J$-invariant part of $S^2\Omega^2M$ to $\Omega^{1,1}M$ is not an isometry. In fact, $|S|^2=2|S\circ J|^2$ (for example $|g|^2=n=2m$ while $|\omega|^2 =m$).  Since the Riemannian curvature $R$ has values in $\Omega^{1,1}M$, the $J$-invariant part of $\Omega^2M$, it can be viewed as a section of the suborbi-bundle of \emph{abstract K\"ahlerian curvature tensors}, $\mathcal{K}M$, defined as the intersection of $\mathcal{R}M$ with $\Omega^{1,1}M \otimes   \Omega^{1,1}M$. However, in general, none of the components of the decomposition \eqref{eq:decomp} is in $\mathcal{K}M$, 
and we have a new decomposition of $R$ inside $\mathcal{K}M$
\begin{equation*}
R=\frac{\scal}{2m(m+1)}(Id_{\vert_{\Omega^{J,+}}} + \omega \otimes \omega) + \frac{1}{m+2}(\{Ric_0,\cdot\}_{\vert_{\Omega^{J,+}}} + \rho_0 \otimes \omega + \omega \otimes \rho_0) + W^{\mathcal{K}},
\end{equation*}
where $\vert_{\Omega^{J,+}}$ is the orthogonal projection of $\Omega^2M$ onto its $J$-invariant part and $W^{\mathcal{K}}$ is the so-called Bochner tensor.

Alternatively, we have (cf. \cite{B}, p. $77$)
\begin{equation}
\label{eq:Bdecomp}
R=\frac{\scal}{2m^2} \omega \otimes \omega + \frac{1}{m} \rho_0 \otimes \omega + \frac{1}{m} \omega \otimes \rho_0 + B,
\end{equation} 
where $B$ decomposes as 
$$
B=\frac{\scal}{2m(m+1)}Id_{\vert_{\Omega^{1,1}M}} + B_0.
$$
Both decompositions \eqref{eq:decomp} and \eqref{eq:Bdecomp} are related by
\begin{align*}
&|B_0|^2  =-\frac{3(m-1)}{m+1}|U|^2 - \frac{m-2}{m}|Z|^2+|W|^2 \\
&|\rho_0|^2  =(m-1)|Z|^2 \\
& \scal^2  = 4m(2m-1) |U|^2.
\end{align*}

In the particular case where $m=2$, 
\eqref{eq:Bdecomp} becomes
\begin{equation}
\label{eq:B2decomp}
R=\frac{\scal}{8} \omega \otimes \omega + \frac{1}{2} \rho_0 \otimes \omega + \frac{1}{2} \omega \otimes \rho_0 + B,
\end{equation} 
with 
\begin{align}
\label{rel:1} & |B_0|^2  =|W|^2 - |U|^2 \\
\label{rel:2} & |\rho_0|^2  = |Z|^2 \\
\label{rel:3} & \scal^2  = 24 |U|^2.
\end{align}


\section{Topological integrals on weighted projective planes}
\label{sec:integrals}
In this section we use localization in equivariant cohomology to obtain expressions for several integrals, topological invariants of the  weighted projective space $\C P^2(N_1, N_2, N_3)$, in terms of the $N_i$'s. These will be used in the proof of Theorem~\ref{thm:1} in \S \ref{sec:mainthm}.

Let us consider $M: =\C P^2(N_1, N_2, N_3)$ and take the Hamiltonian $S^1$-action from Example~\ref{example:1}. Then Proposition~\ref{prop:2.1} becomes  
\begin{equation}
\label{eq:loc2}
\int_M \alpha = \frac{N_3 \, \iota_{F_3}^* \alpha}{N_1\, N_2\,\textsl{x}^2 } + \int_\Sigma  \frac{\iota_\Sigma^* \alpha}{b_\Sigma \, \textsl{u} + \textsl{x}},
\end{equation}
where: $F_3=[0:0:1]$ is the isolated fixed point (at which the $S^1$-moment map $\phi$ takes its maximum value); $\Sigma=\C P^1(N_1,N_2)$ is the orbi-surface on which $\phi$ takes its minimum value; and $b_\Sigma$ is the degree of $\nu_\Sigma$, the normal orbi-bundle of $\Sigma$, i.e., $b_\Sigma\, \textsl{u} = c_1(\nu_\Sigma)$ for a generator $\textsl{u}\in H^2(\Sigma)$ (cf. \S \ref{sec:2}). Here we used the fact that the orbi-weight over a fiber of $\nu_\Sigma$ is equal to $\textsl{x}$, implying that $c_1^{S^1}(\nu_\Sigma)=b_\Sigma\, \textsl{u}+ \textsl{x}$. 

Applying  (\ref{eq:loc2}) to $\alpha = 1$ yields 
$$
0= \frac{N_3}{N_1 \, N_2 \, \textsl{x}^2} +  \int_\Sigma  \frac{1}{b_\Sigma \, \textsl{u} + \textsl{x} } =  \frac{N_3}{N_1 \, N_2 \, \textsl{x}^2} + \frac{1}{ \textsl{x}} \int_\Sigma \sum_{j=0}^\infty (-1)^j\, \left( \frac{b_\Sigma\textsl{u}}{\textsl{x}} \right)^j ,
$$
and we conclude, as expected, that $b_\Sigma = \frac{N_3}{N_1 \, N_2}$. 

Using (\ref{eq:loc2}) with the equivariant symplectic form $\omega^\sharp:=\omega - \phi\, \textsl{x}$ (cf. \cite{AB}) gives
\begin{eqnarray*}
0 = \int_M \omega^\sharp & = &   \frac{N_3 \, \iota_{F_3} ^*\omega^\sharp}{N_1 \, N_2 \, \textsl{x}^2} +  \int_\Sigma  \frac{\iota_\Sigma^* \omega^\sharp}{b_\Sigma \, \textsl{u} + \textsl{x} }
   =   - \frac{N_3}{N_1\, N_2 \textsl{x}}\, y_\M +  \int_\Sigma \frac{\omega - y_\m \, \textsl{x}}{b_\Sigma \, \textsl{u} + \textsl{x} }
\end{eqnarray*}
where $y_\M$ and $y_\m$ denote the values $\phi(F_3)$ and $\phi(\Sigma)$, i.e., the maximum and minimum values of $\phi$. Hence, 
$$
0 = - \frac{N_3}{N_1\, N_2 \textsl{x}}\, y_\M + \int_\Sigma (\frac{\omega}{\textsl{x}} - y_\m)\, \sum_{j=0}^\infty (-1)^j \, \left(\frac{b_\Sigma\, \textsl{u}}{\textsl{x}}\right)^j
$$
and we have
\begin{equation}
\label{eq:area}
\mathrm{area}(\Sigma)=b_\Sigma (y_\M - y_\m)= \frac{N_3}{N_1 \, N_2}(y_\M - y_\m).
\end{equation}

Moreover, using (\ref{eq:loc2}) with $\alpha= (\omega^\sharp)^2$ yields

\begin{eqnarray*}
\int_M (\omega^\sharp)^2 &  = &  \frac{N_3 (\iota_F^*\omega^\sharp)^2}{N_1 \, N_2 \, \textsl{x}^2} +  \int_\Sigma  \frac{(\iota_\Sigma^* \omega^\sharp)^2}{b_\Sigma \, \textsl{u} + \textsl{x} }
  =  \frac{N_3}{N_1\, N_2}\, y_\M^2 +  \int_\Sigma \frac{(\omega - y_\m \, \textsl{x})^2}{b_\Sigma \, \textsl{u} + \textsl{x} }\\
  & = & \frac{N_3}{N_1\, N_2}\, y_\M^2 + \frac{1}{\textsl{x}} \int_\Sigma (-2 y_\m\, \omega\, \textsl{x} + y_\m^2\, \textsl{x}^2)\, \sum_{j=0}^\infty (-1)^j \, \left(\frac{b_\Sigma \, \textsl{u}}{\textsl{x}}\right)^j \\ & = & \frac{N_3}{N_1\, N_2}\, y_\M^2 - 2 y_\m \, \mathrm{area}(\Sigma) - y_\m^2\, b_\Sigma \\
\\ & = & \frac{N_3}{N_1\, N_2}\,(y_\M^2  - 2 \,y_\m\, (y_\M - y_\m) -  y_\m^2) \\ & = &   \frac{N_3}{N_1\, N_2}\,(y_\M -  y_\m)^2.
\end{eqnarray*}
On the other hand, since 
$$
\int_M (\omega^\sharp)^2 = \int_M(\omega - \phi \textsl{x})^2 = \int_M \omega^2 =2 \, \mathrm{vol}(M),
$$
we conclude that 
\begin{equation}
\label{eq:volume}
2\, \mathrm{vol}(M)= \frac{N_3}{N_1 N_2}(y_\M - y_\m)^2.
\end{equation}

If instead we use (\ref{eq:loc2}) with $\alpha= c_1^{S^1}(TM)\wedge \omega^\sharp $, we obtain (cf. \S \ref{sec:2.3})
\begin{eqnarray*}
 \int_M c_1^{S^1}\wedge \omega^\sharp & = &   \frac{N_1 + N_2}{N_1\, N_2}y_\M + \int_\Sigma \frac{(c_1(T\Sigma) + b_\Sigma\, \textsl{u} + \textsl{x})(\omega - \phi \textsl{x})}{b_\Sigma\, \textsl{u} + \textsl{x}} \\
  & = &    \frac{N_1 + N_2}{N_1\, N_2}y_\M + \int_\Sigma \frac{((\frac{1}{N_1}+\frac{1}{N_2})\textsl{u}  + b_\Sigma\, \textsl{u} + \textsl{x})(\omega - y_\m \textsl{x})}{b_\Sigma\, \textsl{u} + \textsl{x}}.
\end{eqnarray*}
Here we used the fact that 
$$\iota_\Sigma^* c_1^{S^1}(TM)= c_1^{S^1}(T\Sigma) + b_\Sigma \, \textsl{u} + \textsl{x}$$ 
and that,  for a general orbi-surface $\Sigma$ with cone singularities,
$$
c_1^{S^1}(T\Sigma) = c_1(T\Sigma) = \chi(\Sigma)\, \textsl{u} = (2-2g+\sum_{i=1}^k (\frac{1}{\alpha_i} - 1))\, \textsl{u},
$$
where the $\alpha_i$'s are the orders of  orbifold structure groups of the singularities of $\Sigma$, $g$ is the genus of the underlying topological surface $|\Sigma |$ and $ \chi(\Sigma)$ is the orbifold Euler characteristic of $\Sigma$ (cf. \cite{FS}). Then we have 
\begin{eqnarray*}
 \lefteqn{\int_M c_1^{S^1}\wedge \omega^\sharp = \frac{N_1 + N_2}{N_1\, N_2}\, y_\M + 
\int_\Sigma \omega -  \frac{N_1 + N_2}{N_1 N_2} \int_\Sigma \frac{y_\m \, \textsl{x} \, \textsl{u}}{b_\Sigma \textsl{u} + \textsl{x}} =} \\ & = &  \frac{N_1 + N_2}{N_1\, N_2}\, y_\M  + \mathrm{area}(\Sigma) -  y_\m \frac{N_1 + N_2}{N_1\, N_2} \int_\Sigma \sum_{j=0}^\infty (-1)^j \left(\frac{b_\Sigma}{\textsl{x}}\right)^j \textsl{u}^{j+1} \\ & = & \frac{N_1 + N_2}{N_1\, N_2}\, (y_\M - y_\m) + \mathrm{area}(\Sigma) = \frac{N_1 + N_2 + N_3}{ N_1 N_2}\, (y_\M - y_\m).
\end{eqnarray*}
On the other hand, since $c_1^{S^1}(TM)=c_1(TM) + f \textsl{x}$ where $c_1$ is the ordinary first Chern class and $f$ is a global function on $M$ (cf. \cite{AB}), we have, for dimensional reasons,
$$
\int_M c_1^{S^1} \wedge \omega^\sharp = \int_M c_1 \wedge \omega,
$$ 
and, using \eqref{eq:volume}, we obtain
\begin{equation}
\label{eq:c1w}
\int_M c_1 \wedge \omega = \frac{N_1 + N_2 + N_3}{ N_1 N_2}\, (y_\M - y_\m) = \frac{N_1 + N_2 + N_3}{ \sqrt{N_1 N_2 N_3}}\,  \sqrt{2 \, \mathrm{vol}(M)}.
\end{equation}

Applying again (\ref{eq:loc2}) now to $\alpha= (c_1^{S^1})^2(TM)$ yields 
\begin{eqnarray*}
 \int_M (c_1^{S^1})^2 & = &  \frac{(N_1 + N_2)^2}{N_1\, N_2\,N_3}+ \int_\Sigma \frac{(c_1(T\Sigma) + b_\Sigma\, \textsl{u} + \textsl{x})^2}{b_\Sigma\, \textsl{u} + \textsl{x}} \\
  & = &    \frac{(N_1 + N_2)^2}{N_1\, N_2\, N_3} + \int_\Sigma \frac{((\frac{1}{N_1}+\frac{1}{N_2})\textsl{u}  + b_\Sigma\, \textsl{u} + \textsl{x})^2}{b_\Sigma\, \textsl{u} + \textsl{x}} \\ & = &  \frac{(N_1 + N_2)^2}{N_1\, N_2\,N_3} + 2\left( \frac{1}{N_1} + \frac{1}{N_2}\right) + \int_\Sigma b_\Sigma \, \textsl{u} + \textsl{x} \\ & = & \frac{(N_1+N_2+N_3)^2}{N_1\,N_2\,N_3}.
\end{eqnarray*}
Moreover, since we also have $\int_M (c_1^{S^1})^2 =\int_M c_1^2$, we conclude that
\begin{equation}
\label{eq:c1c1}
\int_M c_1^2 = \frac{(N_1+N_2+N_3)^2}{N_1\,N_2\,N_3}.
\end{equation}

Finally, using (\ref{eq:loc2}) with $\alpha= c_2^{S^1}(TM)$, we get
\begin{eqnarray*}
 \int_M c_2^{S^1} & = &  \frac{1}{N_3}+ \int_\Sigma \frac{c_1(T\Sigma) \cdot (b_\Sigma\, \textsl{u} + \textsl{x})}{b_\Sigma\, \textsl{u} + \textsl{x}} \\
  & = &    \frac{1}{N_3} + \int_\Sigma \frac{(\frac{1}{N_1}+\frac{1}{N_2})\textsl{u} \cdot (b_\Sigma\, \textsl{u} + \textsl{x})}{b_\Sigma\, \textsl{u} + \textsl{x}} \\ & = &  \frac{1}{N_3} + \frac{1}{N_1} + \frac{1}{N_2}. 
\end{eqnarray*}
Here we used the fact that 
$$\iota_\Sigma^* c_2^{S^1}(TM)= \iota_\Sigma^* e(TM) = e(T\Sigma) \cdot e(\nu_\Sigma)= c_1^{S^1}(T\Sigma) \cdot (b_\Sigma \, \textsl{u} + \textsl{x}).$$
Since, on the other hand,  $\int_M c_2^{S^1} =\int_M c_2$, we obtain
\begin{equation}
\label{eq:c2}
\int_M c_2 = \frac{1}{N_1} + \frac{1}{N_2} + \frac{1}{N_3}.
\end{equation}


\section{Scalar curvature of extremal metrics on weighted projective
  planes}
\label{sec:A}

In this section we follow~\cite{Ab} to give a convenient description
of the scalar curvature of extremal K\"ahler metrics on weighted projective planes, using toric geometry. We compute in particular the
square of its $L^2$-norm on $\C P^2(N_1,N_2,N_3)$, 
which is determined by the weights $N_1, N_2, N_3$. 

\subsection{Weighted and labeled projective spaces}
\label{sec:A2}

The weighted projective space $\C P^{m}(\bN)$
described in Example \ref{example:0} is in fact a compact 
\emph{complex} orbifold. 
Consider the finite group $\Ga_\bN$ defined by
$$
\Ga_\bN = \left(\Z_{\hN_1}\times\cdots\times\Z_{\hN_{m+1}}\right)/
\Z_\hN\,,
$$
where 
\[
\hN_r = \prod_{k=1,k\neq r}^{m+1} N_k \quad, \quad \quad
\hN = \prod_{k=1}^{m+1} N_k
\]
and
\begin{eqnarray*}
\Z_\hN & \hookrightarrow & \Z_{\hN_1}\times\cdots\times\Z_{\hN_{m+1}}\\
\zeta & \mapsto & \left(\zeta^{N_1}\,,\  \ldots\ ,\,\zeta^{N_{m+1}}\right)
\end{eqnarray*}
($\Z_q\equiv \Z/q\Z$ is identified with the group of $q$-th roots
of unity in $\C$). $\Ga_\bN$ acts on $\C P^m (\bN)$ via
\[
[\eta]\cdot [z] = [\eta_1 z_1: \ldots : \eta_{m+1} z_{m+1}]\,,\ 
\mbox{for all}\ [\eta]\in\Ga_\bN\,,\ [z]\in\C P^m (\bN)\,,
\]
and we define the labeled projective space $\C P^m [\bN]$ as the
quotient
\[
\C P^m [\bN] := \C P^m (\bN) / \Ga_\bN\,,
\]
and denote its points by $[[z_1: \ldots : z_{m+1}]] \in \C P^m [\bN]$.

By definition, the induced orbifold structure on $\C P^m [\bN]$ is
such that the quotient map
\[
\pi_\bN : \C P^m (\bN) \to \C P^m [\bN]
\]
is an orbifold covering map. This means that any orbifold geometric
structure on $\C P^m [\bN]$ (e.g., symplectic, complex, K\"ahler, etc)
lifts through $\pi_\bN$ to a $\Ga_\bN$-invariant orbifold geometric
structure on $\C P^m (\bN)$.

The action of $\Ga_\bN$ is free on
\[
\breve{\C P}^m (\bN) := \{ [z_1: \ldots : z_{m+1}]\in\C P^m (\bN)\,:\ 
\text{$z_k \ne 0$ for all $k$} \}\,.
\]
In particular, $\pi_\bN$ has degree $|\Ga_\bN| = (\hN)^{m-1}$ and 
\[
\breve{\C P}^m [\bN] := \{ [[z_1: \ldots : z_{m+1}]]\in\C P^m [\bN]\,:\ 
\text{$z_k \ne 0$ for all $k$} \}
\]
is an open dense smooth subset of $\C P^m [\bN]$. On the other hand,
one can check (see~\cite{Ab}) that the orbifold structure group of any
point in
\[
\C P^m [\bN]_r := \{ [[z_1 : \ldots : z_{m+1}]]\in\C P^m [\bN]\,:\ 
\text{$z_r = 0$ and $z_k \ne 0$ for all $k\ne r$} \}
\]
is isomorphic to $\Z_{\hN_r}$.

Both $\C P^m (\bN)$ and $\C P^m[\bN]$ are examples of compact
K\"ahler toric orbifolds, i.e., K\"ahler orbifolds of real dimension
$n=2m$ equipped with an effective holomorphic and Hamiltonian action of
the standard real $m$-torus $\T^m = \R^m / 2\pi\Z^m$. In fact, 
the standard K\"ahler structure and $\T^{m+1}$-action on $\C^{m+1}$ 
induce a suitable K\"ahler structure and $\T^m$-action on these
quotients (see~\cite{Ab}). For the labeled projective space $\C P^m
[\bN]$ the moment map of the Hamiltonian $\T^m$-action, 
$\phi : \C P^m [\bN] \to (\R^m)^\ast$, has image the simplex
$P_\la^m \subset (\R^m)^\ast$ defined by
\[
P_\la^m = \bigcap_{r=1}^{m+1} \left\{ x \in (\R^m)^\ast :
\ell_r(x):= \langle x,\mu_r\rangle - \la \geq 0 \right\}\,,
\]
where: $\la\in\R^+$, $\mu_r = e_r\,,\
r=1,\ldots,m$, $\mu_{m+1} = - \sum_{j=1}^m e_j$ and 
$(e_1,\ldots,e_m)$ denotes the standard basis of $\R^m$ 
(see Figure 1 for the case $m=2$). The positive real number 
$\la\in\R^+$ parametrizes the cohomology
class of the K\"ahler form on $\C P^m [\bN]$.

\begin{figure}[h!] 
  \centering
  \includegraphics[scale=.7]{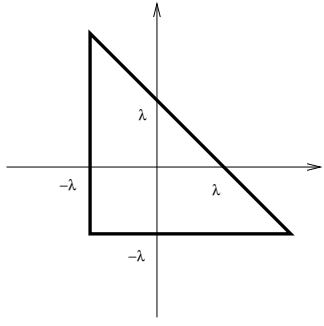}\\
  \caption{The simplex $P_\lambda^2 = \phi(\C P^2 [\bN])$.}
\end{figure}

\subsection{Scalar curvature of extremal metrics}
\label{sec:A3}
 
As shown by Bryant in~\cite{Br}, weighted projective spaces have
Bochner-K\"ahler metrics, i.e., K\"ahler metrics with vanishing Bochner
curvature, which are unique in each K\"ahler class. These metrics are
in particular extremal in the sense of Calabi~\cite{Ca}, i.e. minimize
the $L^2$-norm of the scalar curvature within K\"ahler metrics in a
fixed K\"ahler class (see also \cite{H}). Moreover, they are always
toric and turn out to have a very simple description in terms of
explicit data on the simplex $P_\la^m$ (see~\cite{Ab}). 

In particular, the scalar curvature $\tau_\la$ of these extremal 
Bochner-K\"ahler metrics $\om_\la$, being a $\T^m$-invariant function
on $\C P^m (\bN)$ and $\C P^m [\bN]$, descends to a function on the 
simplex $P_\la^m$ which is given by
\begin{equation}
\tau_\la (x) = \frac{1}{\la}\left( \frac{n}{m+1} 
\left( \sum_{r=1}^{m+1} \frac{1}{\hN_r}\right) +
\frac{2(m+2)}{m+1} \sum_{j=1}^m (\frac{1}{\hN_{m+1}} - \frac{1}{\hN_j})
\frac{x_j}{\la} \right)\,.
\end{equation}

Since the orbifold covering map $\pi_\bN : \C P^m (\bN) \to \C P^m
[\bN]$ has degree $(\hN)^{m-1}$ and the push-forward by the moment map 
$\phi$ of the volume form $\mu_\la := \om_\la^m / m!$ on $\C P^m
[\bN]$ is $(2\pi)^m dx$ on $P^m_\la$, we have that
\begin{eqnarray*}
\frac{1}{4\pi} \int_{\C P^m (\bN)} \tau_\la\,\mu_\la & = &
\frac{1}{4\pi} (\hN)^{m-1} \int_{\C P^m [\bN]} \tau_\la\,\mu_\la \\
& = & \frac{1}{4\pi} (\hN)^{m-1} (2\pi)^m \int_{P_\la^m}
\tau_\la(x)\,dx\,.
\end{eqnarray*}
Using the fact that
\[
\int_{P_\la^m} x_j\,dx = 0\,,\ \forall\, \la\in\R^+, j=1,\ldots,m\,,
\]
we then have
\begin{eqnarray*}
\frac{1}{4\pi} \int_{\C P^m (\bN)} \tau_\la\,\mu_\la & = &
(2\pi)^{m-1} (\hN)^{m-1} \frac{m}{\la(m+1)} 
\left( \sum_{r=1}^{m+1} \frac{1}{\hN_r}\right) \vol(P_\la^m) \\
& = & (2\pi)^{m-1} (\hN)^{m-2} \frac{m}{m+1}
\left( \sum_{r=1}^{m+1} N_r\right) \frac{\vol(P_\la^m)}{\la} \\
& = & (2\pi)^{m-1} (\hN)^{m-2} \frac{m}{m+1}
\left( \sum_{r=1}^{m+1} N_r\right) \la^{m-1} \vol(P_1^m) \\
& = & \frac{\vol(P_1^m)^{\frac{1}{m}}}{\hN^{\frac{1}{m}}} \frac{m}{m+1} 
\left( \sum_{r=1}^{m+1} N_r\right) \left((2\pi\la)^{m}
\hN^{m-1} \vol(P_1^m)\right)^{\frac{m-1}{m}} \\
& = & \frac{\vol(P_1^m)^{\frac{1}{m}}}{\hN^{\frac{1}{m}}} \frac{m}{m+1} 
\left( \sum_{r=1}^{m+1} N_r\right) 
\left( \vol_\la (\C P^m (\bN))\right)^{\frac{m-1}{m}}\,.
\end{eqnarray*}
When $m=2$ we have that $\vol(P_1^2) = 9/2$ and so
\[
\int_{\C P^2(\bN)} c_1 \wedge \om_\la = \frac{1}{4\pi}
\int_{\C P^2(\bN)} \tau_\la\,\mu_\la =  
\frac{N_1+N_2+N_3}{\sqrt{N_1 N_2 N_3}} 
\sqrt{2 \vol_\la (\C P^2(\bN))}\,,
\]
which agrees with~(\ref{eq:c1w}) as expected. Note that our
normalization of the scalar curvature is the usual one in Riemannian 
geometry, which is twice the one often used in K\"ahler geometry.

When $m=2$ the integral of the square of the scalar curvature is scale
invariant, and so
\[
\int_{\C P^2(\bN)} (\tau_\la)^2\,\mu_\la =
\int_{\C P^2(\bN)} (\tau_1)^2\,\mu_1\,.
\]
The scalar curvature can be written in this case as
\[
\tau_1(x) = \frac{4}{3 N_1 N_2 N_3}
\left((N_1 + N_2 + N_3) + 2N_3 (x_1+x_2) - 2N_1 x_1 - 2N_2 x_2 \right)
\,. 
\]
Using the fact that
\[
\int_{P_1^2} (x_i)^2 \,dx = \frac{9}{4} \quad\text{and}\quad
\int_{P_1^2} x_i x_j \,dx = -\frac{9}{8}\,,\ i\ne j\,,
\]
one gets
\begin{eqnarray*}
\int_{\C P^2(\bN)} (\tau_\la)^2\,\mu_\la & = &
(N_1 N_2 N_3) (2\pi)^2 \int_{P_1^2} (\tau_1(x))^2\,dx\\
& = & 96 \pi^2 \frac{N_1^2 + N_2^2 + N_3^2}{N_1 N_2 N_3}\,.
\end{eqnarray*}


\section{Heat Invariants}\label{sec:heat}

To study the relationship between the geometry and the Laplace spectrum of orbifolds with isolated singularities, we will consider the asymptotic expansion of the heat trace $K(x,x,t)$ as $t \rightarrow 0^+$.  The so-called \emph{heat invariants} appearing in this expansion involve geometric quantities such as the dimension, the volume, and the curvature.  For good orbifolds (i.e., those arising as global quotients of manifolds), Donnelly \cite{Don76} proved the existence of the heat kernel and constructed the asymptotic expansion for the heat trace.  This work was extended to general orbifolds in \cite{DGGWZ}, where the expressions obtained also clarify the contributions of the various pieces of the singular set.  
While \cite{Don76} and \cite{DGGWZ} treat the heat trace asymptotics for functions, we will also require the asymptotics for 1-forms.  We begin with the function case, extracting the necessary background from \cite[\S 4]{DGGWZ} and simplifying the expressions when possible to reflect the fact that we are only interested in isolated singularities. 

\begin{definition}\label{def.bk}\cite{Don76}
Let $h$ be an isometry of a Riemannian manifold $X$
and let $\Omega(h)$ denote the set of (isolated) fixed points of $h$.  For $x \in \Omega(h)$,
define $A_h(x):=h_*:T_x(X) \to T_x(X)$; note that $A_h(x)$ is non-singular.  
Set $$B_h(x)=(I-A_h(x))^{-1}.$$  
\end{definition}

\begin{proposition}\label{prop.don}\cite{Don76}  Let $(X,g)$ be a closed Riemannian manifold,
let $K(t,x,y)$ be the heat kernel of $X$, and let $h$ be a nontrivial isometry
of $X$ with isolated fixed points. 
Then $\int_X\,K(t,x,h(x)) \, \mu_g (x)$ has
an asymptotic expansion as $t\rightarrow 0^+$ of the form
$$\sum_{x\in\Omega(h)}
\sum_{k=0}^\infty\,t^k \tilde{b}_k(h).$$   
For $x \in \Omega(h)$,   $\tilde{b}_k(h)$ depends only  on the germs of
$h$ and of the Riemannian metric of $X$ at $x$.
 \end{proposition}

Donnelly gives explicit formulae for $\tilde{b}_0$ and $\tilde{b}_1$ \cite[Thm. 5.1]{Don76} as follows.  All indices run from $1$ to $n$, where $n$ is the dimension of the Riemannian manifold $X$.  
At each point $x\in \Omega(h)$, choose an orthonormal basis
$\{e_1,\dots,e_n\}$ of $T_x(X)$.  The sign convention on the curvature tensor 
$R$ of $X$ is chosen so that $R_{abab}$  is the sectional curvature of the 2-plane 
spanned by $e_a$ and $e_b$.  Set
$$\scal=\sum_{a,b=1}^n\,R_{abab}$$
and $$Ric_{ab}=\sum_{c=1}^n\,R_{acbc}.$$
Thus $\scal$ is the scalar curvature and $Ric$ the Ricci tensor of $X$.  Then
\begin{equation}\label{eq.b0}\tilde{b}_0(h)=|det(B_h(x))|\end{equation}
and, summing over repeated indices, 
\begin{equation}\label{eq.b1}\begin{aligned}
\tilde{b}_1(h)=|det(B_h(x))|&\left( \frac{1}{6}\scal+\frac{1}{6}Ric_{kk}+\frac{1}{3}R_{iks
j}B_{ki}B_{js} \right. \\
 & \left. +\frac{1}{3}R_{iktj}B_{kt}B_{ji}-R_{kaja}B_{ks}B_{js} \right).
\end{aligned}\end{equation}

We want to see what the analogous expressions are for orbifolds which are not necessarily global quotients.  We again restrict to the case of orbifolds which have only isolated singularities and present the appropriately simplified expressions.  Let $M$ be such an orbifold, and choose a singularity $x \in M$.  Let $(U,V,\Gamma)$ be an orbifold chart for a neighborhood of $x$, and let $\gamma \in
\Gamma$.  Define
$$
b_k(\gamma) = \tilde{b}_k(\sigma(\gamma)),
$$
where $\sigma$ is an isomorphism from $\Gamma$ to the isotropy group of a point in the preimage of $x$ under the homeomorphism given by the orbifold chart.  That is, we use the charts to calculate the value of $b_k$ locally.  Set
$$I_{\gamma}:= \sum_{k=0}^\infty\,t^k \,b_k(\gamma),$$
and
$$I_x:=\sum_{\gamma \in \Gamma}\,I_{\gamma}.$$
Also set
$$I_0:=(4\pi t)^{-dim(M)/2}\sum_{k=0}^\infty\,a_k(M)t^k$$
where the $a_k(M)$ (which we will usually write simply as $a_k$) are the
familiar heat invariants.  In particular, $a_0=vol(M)$,
$a_1=\frac{1}{6}\int_M\scal(x) \, \mu_g$, etc.   Note that if $M$ is finitely
covered by a Riemannian manifold $X$, say $M=G\backslash X$, then
$a_k(M)=\frac{1}{|G|}a_k(X)$.

We can now give an expression for the asymptotic expansion of the heat trace of an orbifold $M$ with isolated singularities (see \cite{DGGWZ} for more details).

\begin{Theorem}\label{thm.trace}\cite{DGGWZ} 
Let $M$ be a Riemannian orbifold and let
$\lambda_1\leq \lambda_2\leq\dots$ be the spectrum of the associated Laplacian
acting on smooth functions on $M$.  The heat trace
$\sum_{j=1}^\infty\,e^{-\lambda_jt}$ of $M$ is asymptotic as $t \rightarrow 0^+$ to 
\begin{equation}\label{eq.as}
I_0+\sum_{x \in \Omega}\,\frac{I_x}{|\Gamma|}
\end{equation}
where $\Omega$ is the set of isolated singularities of $M$.  This asymptotic expansion is of
the form
\begin{equation}\label{eq.pow}
(4\pi t)^{-dim(M)/2}\sum_{j=0}^\infty\,c_jt^{\frac{j}{2}}.
\end{equation}
\end{Theorem}

Using this expression, we calculate the first few terms in the asymptotic expansion of the heat trace of a 4-dimensional Riemannian orbifold $M$ with isolated singularities.  Let $x \in M$ be a cone point of order $N$.  For any chart $(U,V,\Gamma)$ about $x$, let $\gamma$ generate $\Gamma$; for $j=1,\dots,N-1$ we have (cf. \S \ref{sec:2.3})
$$A_{\gamma^j}=\gamma^j_*=
\begin{bmatrix} 
\cos(\frac{2j\pi}{N}) & -\sin(\frac{2j\pi}{N}) & 0 & 0 \\
\sin(\frac{2j\pi}{N}) & \cos(\frac{2j\pi}{N}) & 0 & 0 \\
0 & 0 & \cos(\frac{2mj\pi}{N}) & -\sin(\frac{2mj\pi}{N}) \\
0 & 0 & \sin(\frac{2mj\pi}{N}) & \cos(\frac{2mj\pi}{N})
\end{bmatrix},$$
for some integer $1 \leq m < N$ with $(m,N)=1$.
Thus
$$b_0(\gamma^j)=|det((I-A_{\gamma^j})^{-1})|= \frac{1}{16 \sin^2(\frac{j \pi}{N})\sin^2(\frac{mj\pi}{N})}.$$
Hence
$$
I_x = \sum_{j=1}^{N-1} \frac{1}{16 \sin^2(\frac{j \pi}{N})\sin^2(\frac{mj\pi}{N})}+ O(t).
$$

Now consider our weighted projective space $M:=\mathbb{C}P^2(N_1,N_2,N_3)$.  Note that 
\begin{eqnarray*}
I_0 & = & (4\pi t)^{-2}\sum_{k=0}^\infty\,a_k(M)t^k \\
 & = & \frac{1}{16 \pi^2} (a_0 t^{-2} + a_1 t^{-1} + a_2 + a_3 t + \cdots).
\end{eqnarray*}
By Theorem \ref{thm.trace}, we see that the heat trace of $M$ is asymptotic as $t \rightarrow 0^+$ to
$$
\frac{1}{16 \pi^2} (a_0 t^{-2} + a_1 t^{-1} + a_2) + T + O(t),
$$
where 
\begin{eqnarray*}
T & =  & \frac{1}{N_1}\sum_{j=1}^{N_1-1} \frac{1}{16 \sin^2(\frac{N_2j \pi}{N_1})\sin^2(\frac{N_3j\pi}{N_1})} + \frac{1}{N_2}\sum_{j=1}^{N_2-1} \frac{1}{16 \sin^2(\frac{N_1j \pi}{N_2})\sin^2(\frac{N_3j\pi}{N_2})} \\
& & + \frac{1}{N_3}\sum_{j=1}^{N_3-1} \frac{1}{16 \sin^2(\frac{N_1j \pi}{N_3})\sin^2(\frac{N_2j\pi}{N_3})}
\end{eqnarray*}
and we have replaced the factor $m$ by the appropriate ratios of the $N_i$'s (see Example 2).
Hence the coefficient of the term of degree -2 is $\frac{a_0}{16 \pi^2}$, the coefficient of the term of degree -1 is $\frac{a_1}{16 \pi^2}$, and that of the term of degree 0 is $\frac{a_2}{16\pi^2} + T$.  This means that the spectrum determines
\begin{eqnarray*}
a_0 & = & vol(M) \\
a_1 & = & \frac{1}{6}\int_M\scal(x) \, \mu_g 
\end{eqnarray*}
and 
$$
\frac{1}{360*16 \pi^2}\int_M (2 ||R||^2 - 2|Ric|^2 + 5\scal^2) \, \mu_g + T.
$$

\begin{remark}
\label{rmk:1}
Here, the norm $||.||$ is obtained by contracting the tensor with itself, while the norm $|.|$ used in \cite{B} and in \S \ref{sec:2.4} is the usual norm for antisymmetric products of symmetric tensors. Therefore, $||T \wedge S|| = 2 |T \wedge S|$ for any pair of rank-$2$ symmetric tensors and $||T||=|T|$. In particular, $||R|| = 2 |R|$, while $||Ric||=|Ric|$. 
\end{remark}

To complete our study of the heat trace asymptotics, we briefly discuss the case of forms.  As Donnelly notes in \cite{Don79}, the results therein ``generalize easily to the Laplacian with coefficients in a bundle.''  Indeed, by examining \cite{Don76,DonPa}, we see that the appropriate coefficient $k(p)$ on the singular part of the heat expansion for forms is $\binom{n}{p}$, where $n$ is the dimension of the orbifold $M$ and $p$ is the level of form.  On the smooth part of the expansion, Theorem 3.7.1 of \cite{Gil04} gives a coefficient of $k(p)= \binom{n}{p}$ on $a_0$, a coefficient of $k_0(p) = \binom{n}{p} - 6 \binom{n-2}{p-1}$ on $a_1$, and coefficients on the curvature terms appearing in $a_2$ of $k_1(p) = 2 \binom{n}{p} - 30 \binom{n-2}{p-1} + 180 \binom{n-4}{p-2}$, $k_2(p) = -2\binom{n}{p} + 180 \binom{n-2}{p-1} - 720 \binom{n-4}{p-2}$, $k_3(p) = 5 \binom{n}{p} - 60 \binom{n-2}{p-1} + 180\binom{n-4}{p-2}$ on $||R||^2, |Ric|^2,$ and $\scal^2$, respectively (i.e., $a_2 = \frac{1}{360}(4 \pi)^{-n/2} \int_M \{k_1(p)||R||^2 + k_2(p)|Ric|^2 + k_3(p)\scal^2 \} \, \mu_g$).  
For $n=4$, the following table gives the values for these coefficients for $p=0$ and $p=1$; we will use these values in computations in \S \ref{sec:mainthm}.

\begin{table}[h!]\label{table:1}
\begin{tabular}{|c|c|c|c|c|c|} 
\hline  & $k(p)$ &  $k_0(p)$ & $k_1(p)$ & $k_2(p)$ & $k_3(p)$ \\ \hline \hline
$p=0$ & $1$ & $1$ & $2$ & $-2$ & $5$ \\ \hline
$p=1$ & $4$ & $-2$ & $-22$ & $172$ & $-40$ \\ \hline
\end{tabular}
\end{table}


\section{Hearing the weights and extremal metrics}
\label{sec:mainthm}

Using the background and tools developed in the preceding sections, we now present and prove several results related to spectral determination of the weights of four-dimensional weighted projective spaces.

\subsection{Proof of Theorem \ref{thm:1}}
Let $R$ be the full curvature of $M$, let $Ric$ be the full Ricci curvature, and let $\scal$ be the scalar curvature.
We know from \S \ref{sec:heat} that the spectrum of the form-valued Laplacian determines
\begin{align}
\mathfrak{a}_0 (\Delta_p) & = \frac{1}{16 \pi^2} \, k(p) \, \mathrm{Vol}(M)\\
\mathfrak{a}_1 (\Delta_p) & = \frac{1}{16 \pi^2} \, \frac{1}{6} \, k_0(p) \int_M \scal \, \mu_g \\
\label{eqn:a_2} \mathfrak{a}_2 (\Delta_p) & = 
\frac{1}{16 \pi^2}\, \frac{1}{360}  \int_M \left(k_1(p) ||R||^2 + k_2(p) |Ric|^2 + k_3(p) \scal^2\right) \, \mu_g + \,\, k(p) \, T,
\end{align}
where $T$ is the trigonometric sum
\begin{align*}
T := & \frac{1}{16} \left( \frac{1}{N_1} \, \sum_{j=1}^{N_1-1} \frac{1}{\sin^2(\frac{N_2 j \pi}{N_1})\sin^2(\frac{N_3 j \pi}{N_1})} + 
\frac{1}{N_2}\, \sum_{j=1}^{N_2-1} \frac{1}{\sin^2(\frac{N_1 j \pi}{N_2})\sin^2(\frac{N_3 j \pi}{N_2})}+  \right. \\
& \left. +  \frac{1}{N_3}\, \sum_{j=1}^{N_3-1} \frac{1}{\sin^2(\frac{N_1 j \pi}{N_3})\sin^2(\frac{N_2 j \pi}{N_3})} \right).
\end{align*}
From decompositions \eqref{eq:decomp} and \eqref{eq:B2decomp} and relations \eqref{rel:1} through \eqref{rel:3} we have (see also Remark~\ref{rmk:1})
\begin{eqnarray*}
\lefteqn{k_1(p) ||R||^2 + k_2(p) |Ric|^2 + k_3(p) \scal^2= }\\ & = & 4 k_1(p) (|U|^2+|Z|^2+|W|^2) + k_2(p) |Ric|^2 + k_3(p) \scal^2 \\
 & = & 4 k_1(p)(|U|^2 + |\rho_0|^2 + |W|^2) + k_2(p)(|Ric_0|^2+\frac{\tau^2}{4})+ k_3(p)\scal^2 \\ & = & \left(\frac{k_1(p)}{6} +\frac{k_2(p)}{4}+k_3(p)\right)  \scal^2 + 2\left(2k_1(p)+k_2(p)\right) |\rho_0|^2 + 4k_1(p)\left(|B_0|^2+|U|^2\right) \\
 & = & \left(\frac{k_1(p)}{3} +\frac{k_2(p)}{4}+k_3(p)\right) \scal^2 + 2\left(2k_1(p)+k_2(p)\right) |\rho_0|^2 + 4k_1(p)|B_0|^2,
\end{eqnarray*}
where we used that $Ric = \frac{\text{tr}\, Ric}{2m}\,g + Ric_0 =\frac{\scal}{4}\,g + Ric_0$ (implying that $|Ric|^2=\frac{\scal^2}{4} + |Ric_0|^2$), and that $2|\rho_0|^2=|Ric_0|^2$ .
Moreover, we know that on a complex orbifold of (real) dimension $4$, 
\begin{equation}
4\pi^2 \int_M c_1^2 =  \int_M (\frac{\scal^2}{8} - |\rho_0|^2 )  \, \, \mu_g 
\end{equation}
and 
\begin{equation}
8 \pi^2 \int_M c_2 = \int_M (\frac{\scal^2}{12} - |\rho_0|^2 + |B_0|^2)  \, \, \mu_g 
\end{equation}
(cf. \cite[p. 80]{B} for the proof in the manifold case), and so (\ref{eqn:a_2}) becomes 
\begin{eqnarray*}
\lefteqn{\mathfrak{a}_2(\Delta_p) = \frac{1}{360 * 16 \pi^2} \left(k_1(p)+\frac{1}{2}\, k_2(p) +k_3(p)\right) \, \int_M \scal^2 \,\,+ } \\  & + & \frac{1}{180}\, k_1(p) \, \int_M  c_2 \,\, - \, \, \frac{1}{720}\left(4 k_1(p)+k_2(p)\right)\, \int_M  c_1^2 \,\, + \,\, k(p)\, T.
\end{eqnarray*}
Using  the values of $k_i(p)$ for $p=0,1$  (cf. Table~\ref{table:1}), we obtain
\begin{equation}\label{eqn:matrix}
\left[ \begin{array}{rrr} \frac{1}{960 \pi^2} & \frac{1}{90} & 1 \\ &  &  \\ \frac{1}{240 \pi^2} & -\frac{11}{90} & 4 \end{array}\right] \left[\begin{array}{c} \int_M \scal^2 \\  \\ \int_M c_2 \\ \\ T \end{array} \right] = \left[\begin{array}{c} \mathfrak{a}_2(\Delta_0) + \frac{1}{120} \int_M c_1^2\\  \\ \mathfrak{a}_2(\Delta_1) +\frac{7}{60} \int_M c_1^2 \end{array} \right],
\end{equation}
implying that 
\begin{equation}
\label{eq:c_2}
\int_M c_2 = - 6 \left(\mathfrak{a}_2(\Delta_1)-4\mathfrak{a}_2(\Delta_0) + \frac{1}{12} \int_M c_1^2 \right).
\end{equation}

On the other hand, the integral of the scalar curvature is a topological invariant, depending only on the K\"ahler class represented by $\omega$ and on the first Chern class in the following way,
\begin{equation*}
\int_M \scal \, \mu_g = 2\pi \int_M c_1 \wedge \omega.
\end{equation*}  
Hence, using for instance $\mathfrak{a}_0(\Delta_0)$ and $\mathfrak{a}_1(\Delta_0)$, we have by \eqref{eq:c1w}
\begin{equation}
\label{eq:a1}
\mathfrak{a}_1(\Delta_0)=  \frac{1}{24 \pi} \sqrt{2 \pi \mathfrak{a}_0(\Delta_0)\,b},
\end{equation}
where  
\begin{equation}
\label{eq:b}
b:= \frac{(N_1+N_2+N_3)^2}{N_1 N_2 N_3} = \int_M c_1^2 
\end{equation}
(cf. \eqref{eq:c1c1}). Therefore, we conclude that  we can hear $b$ from the heat invariants and then, using \eqref{eq:c_2}, we can also hear 
\begin{equation}
\label{eq:c}
c: = \int_M c_2  = \frac{ N_1N_2+N_1 N_3 + N_2N_3}{N_1 N_2 N_3}
\end{equation}
(see \eqref{eq:c2}), as well as
\begin{equation}
\label{eq:a}
d:= \frac{N_1^2 + N_2^2 +N_3^2}{N_1 N_2 N_3} = b - 2c.
\end{equation}

We will now see that these three numbers $b$, $c$, $d$ determine $N_1$, $N_2$ and $N_3$.
First, we note that in $c$, the numerator $N_1 N_2 +  N_1 N_3 + N_2 N_3$ is relatively prime with the denominator
$N_1 N_2 N_3$. Indeed, if these two integers had a common divisor $l$, this would have to divide one
of the $N_i$'s; let us assume without loss of generality that  $l$ divided $N_1$; then $l$ would necessarily divide $N_2 N_3$ (since $l$ would also be a
divisor of $N_1 N_2 + N_1 N_3$) which is impossible since the $N_i$'s are pairwise
relatively prime. Knowing this, we conclude that we can also hear $s:=N_1 N_2 N_3$, the  smallest integer that multiplied by $c$ produces an
integer, and then, multiplying $b$, $c$  and $d$  by this number, we hear the integers:
\begin{align}
\label{eq:p}  p & = N_1 + N_2 + N_3\\
\label{eq:q}  q & = N_1^2 + N_2^2 + N_3^2 \\
\label{eq:r}  r & = N_1 N_2 +  N_1 N_3 + N_2 N_3.
\end{align}
Writing $u=N_1+N_2$ and $v=N_1\,N_2$ we obtain
\begin{align*}
 s & = v(p-u) \\
 r & = u(p-u) + v, 
\end{align*}
implying that
\begin{equation}
\label{eq:3poly}
u^3-2\, p\,u^2 + (p^2 + r)\,u +\,(s-p\,r)=0.
\end{equation}
This equation determines $N_1$, $N_2$ and $N_3$ uniquely up to permutation. Indeed equation  \eqref{eq:3poly} has at most three solutions (three possible values for $N_1+N_2$) which, by \eqref{eq:p} give us the possible values of $N_3$; moreover, we can determine  $N_1\,N_2$ and $N_1^2+N_2^2$ using \eqref{eq:r} and \eqref{eq:q} respectively; then, using these, we can compute $(N_1-N_2)^2$ and consequently  obtain all possibilities for $N_1-N_2$ which, combined with the results for $N_1+N_2$, give us the possible values of $N_1$ and $N_2$.

Note that equation \eqref{eq:3poly} has a unique solution exactly when $p^2-3r=0$. Indeed, for $f(u)=u^3-2\,p\,u^2 + (p^2 + r)\,u +(s-p\,r)$ to have only one zero we need $f^\prime(u)=3\,u^2 - 4\,p\,u + p^2+r$ to have at most one zero and that occurs when the discriminant $p^2-3\,r= 0$. Since  
$$
p^2-3r = N_1^2-N_1(N_2+N_3)+N_2^2+N_3^2-N_2\,N_3, 
$$
this condition is attained when 
$$
N_1= \frac{1}{2}(N_2+N_3 \pm \sqrt{-3\,(N_2-N_3)^2}\,),
$$
implying that $N_1=N_2=N_3$ and then, since they are pairwise relatively prime we have $N_1=N_2=N_3=1$. We conclude that \eqref{eq:3poly} has a unique solution if and only if $M$ is a smooth manifold, thus finishing the proof of Theorem \ref{thm:1}.

\begin{remark}
\label{rmk:3}
Using the methods of \S \ref{sec:integrals} with higher-dimensional weighted
projective spaces $M=\C P^m({\bf N})$ we can obtain the values of  the topological
integrals
\begin{align*}
\int_M c_2 \wedge \omega^{m-2} & = \left(m!
\text{Vol}(M)\right)^{\frac{m-2}{m}}\frac{\sum_{i,j=1}^{m+1} N_iN_j}{(N_1\cdots
N_{m+1})^\frac{2}{m}} \ \ \ ,\\
\int_M c_1^2 \wedge \omega^{m-2} & = \left(m!
\text{Vol}(M)\right)^{\frac{m-2}{m}} \frac{(N_1 + \cdots + N_{m+1})^2}{(N_1\cdots
N_{m+1})^\frac{2}{m}}.
 \end{align*}

Moreover, using $\mathfrak{a}_i(\Delta_p)$ for $i,p=0,1,2$, the curvature decompositions of
\S \ref{sec:2.4} and the expressions
\begin{align*}
\frac{4\pi^2}{(m-2)!}\int_M c_1^2 \wedge \omega^{m-2} & = \int_M \left( \frac{m-1}{4m}
\tau^2 - |\rho_0|^2 \right) \mu_g \\
\frac{8\pi^2}{(m-2)!}\int_M c_2 \wedge \omega^{m-2} & = \int_M
\left( \frac{m-1}{4(m+1)}  \tau^2 - \frac{2(m-1)}{m}|\rho_0|^2 + |B_0|^2 \right) \mu_g
\end{align*}
found for instance in \cite{B}, we conclude that we can hear
\begin{align*}
\int_M \tau^2, & & \int_M c_2 \wedge \omega^{m-2} & & \text{and} & & \int_M c_1^2\wedge
\omega^{m-2}, \end{align*}
implying that we can hear
\begin{align*}
\frac{\sum_{i,j=1}^{m+1} N_i N_j}{(N_1\cdots N_{m+1})^{\frac{2}{m}}} & &
\text{and} & & \frac{(N_1+\cdots + N_{m+1})^{m}}{N_1\cdots N_{m+1}}.
\end{align*}
To obtain more information on the weights from the Laplace spectra, one could impose
restrictions on the metric and use higher-order terms of the asymptotic expansion
of the heat trace.
 \end{remark}

\subsection{Hearing weights from the function spectrum\label{sec:01}}
As mentioned in \S \ref{sec:1}, we are unable to extract the weights using only the information provided by the first few heat invariants for the $0$-form spectrum.  The recourse to $1$-forms was necessitated by the lack of a closed formula for the value of the trigonometric sum $T$ in terms of the $N_i$'s.  
Note that, since
\begin{equation*}
\frac{1}{\sin^2 \alpha \, \sin^2 \beta}=(1+ \cot^2 \alpha)(1+\cot^2 \beta) = 1 + \cot^2 \alpha + \cot^2 \beta + \cot^2 \alpha \cot^2 \beta 
\end{equation*}
and 
\begin{equation*}
\sum_{j=1}^{N-1} \cot^2 \left(\frac{j \pi}{N}\right) = \frac{1}{3}(N-1)(N-3) 
\end{equation*}
(cf. \cite{B-Y}), we have
\begin{tiny}
\begin{eqnarray} \label{eq:T} 
& T & = \frac{1}{16} \left( \frac{1}{N_1} + \frac{1}{N_2} + \frac{1}{N_3} + \frac{(N_1-1)(N_1-2)}{3 N_1}  + \frac{(N_2-1)(N_2-2)}{3 N_2}+ \frac{(N_3-1)(N_3-2)}{3 N_3} \right. \\ \nonumber 
& + &  \left. \frac{1}{N_1}d(N_1;N_2,N_2,N_3,N_3) + \frac{1}{N_2}d(N_2;N_1,N_1,N_3,N_3) + \frac{1}{N_3}d(N_3;N_1,N_1,N_2,N_2) \right) \\\nonumber
 & = & \frac{1}{16} \left( -3 + 3(\frac{1}{N_1} + \frac{1}{N_2} + \frac{1}{N_3}) + \frac{1}{3}(N_1+N_2+N_3) +   \frac{1}{N_1}d(N_1;N_2,N_2,N_3,N_3) +\right. \\ \nonumber& + & \left. \frac{1}{N_2}d(N_2;N_1,N_1,N_3,N_3) + \frac{1}{N_3}d(N_3;N_1,N_1,N_2,N_2) \right),
\end{eqnarray}
\end{tiny}
\hspace{-.17cm} 
where, for positive integers $p_i$, $d(p_0;p_1,p_2,p_3,p_4)$ is the higher-dimensional Dedekind sum 
\begin{equation}
d(p_0;p_1,p_2,p_3,p_4) := \sum_{j=1}^{p_0 - 1} \left( \prod_{i=1}^4 \cot \left( \frac{j \pi p_i}{p_0} \right) \right)
\end{equation}
(see for example \cite{Z}). Nevertheless, the denominators of each Dedekind sum in \eqref{eq:T} are known (see \cite{Z}) and allow us, in most cases, to determine the product $N_1\,N_2\, N_3$. Then, using the expression for $b$ given by \eqref{eq:b}, we are able to determine the sum $N_1+N_2+N_3$. However, these two values alone are not enough to determine the $N_i$'s uniquely.

If we assume the weights $N_1$, $N_2$ and $N_3$ to be prime then we can determine the orders of the singularities using only $\mathfrak{a}_0$ and $\mathfrak{a}_1$ for functions.

\begin{Theorem}
\label{thm:2}
Let $M:=\C P^2(N_1,N_2,N_3)$ be a four-dimensional weighted projective space with isolated singularities.  Assume the weights $N_1$, $N_2$ and $N_3$ are all prime.  Then the spectrum of the Laplacian acting on functions on $M$ determines the weights. 
\end{Theorem}

\begin{proof}
From the values of $\mathfrak{a}_0$ and $\mathfrak{a}_1$ for functions we can determine the value of $b=\frac{(N_1+N_2+N_3)^2}{N_1\, N_2 \, N_3}$ in \eqref{eq:b} using equation \eqref{eq:a1}. 
Note that we can hear the smallest integer which when multiplied by $b$ produces an integer.  Call this smallest integer $D$ and note that it is a product of at most three prime numbers.

If $D$ is a product of three primes then they are necessarily equal to $N_1$, $N_2$ and $N_3$.

If $D$ is a product of two primes, we immediately know the values of the orders of two of the singularities, say $N_1$ and $N_2$. To obtain the value of the third we first multiply $b$ by the product $N_1\, N_2$ to determine the value of 
$$
l:= \frac{(N_1+N_2+N_3)^2}{N_3}\in \Z;
$$
then $N_3$ necessarily divides $(N_1+N_2+N_3)^2$ and consequently divides $N_1+N_2+N_3$.  Thus
$$
l=N_3 \, p_1^{2m_1}\cdots p_k^{2m_k}
$$
for some prime numbers $p_1, \ldots, p_k$ and some positive integers $m_1,\ldots, m_k$; we conclude that we can hear  $(p_1^{m_1}\cdots p_k^{m_k})^2$ (the greatest perfect square that divides $l$) and hence $N_3$.

Similarly, if $D$ is a product of just one prime, say $N_1$, then we know the value of the order of one singularity. To retrieve the values of the other two we multiply $b$ by $N_1$ to obtain  
$$
l:= \frac{(N_1+N_2+N_3)^2}{N_2\, N_3}\in \Z.
$$
Since $(N_2,N_3)=1$ and $l= N_2 \, N_3 \, (p_1^{m_1}\cdots p_k^{m_k})^2$ for some prime numbers $p_1, \ldots, p_k$ and some positive integers $m_1,\ldots, m_k$, we can hear the value of the product $N_2 \, N_3$ and consequently of $N_2$ and $N_3$.

Finally, if $D=1$ (i.e., if $b$ is an integer) then $(N_1+N_2+N_3)^2$ is the product of $N_1 \, N_2 \, N_3$ by a perfect square and we can again determine $N_1 N_2 N_3$ and from that $N_1$, $N_2$ and $N_3$.
\end{proof}

If we remove the restriction that $N_1$, $N_2$ and $N_3$ be prime, we can still determine the orders of the singularities from $\mathfrak{a}_0$ and $\mathfrak{a}_1$ for functions if we fix the Euler characteristic of the orbifold. Indeed, if we are given $\chi(M):=\frac{1}{N_1} + \frac{1}{N_2} + \frac{1}{N_3}$, we know the value of  $c$ in \eqref{eq:c} and so, following the argument in the proof of Theorem~\ref{thm:1}, we are able to determine $s=N_1 \, N_2 \, N_3$. Assuming we know $\mathfrak{a}_0$ and $\mathfrak{a}_1$ we can again determine the value of $b$ in \eqref{eq:b}, using equation \eqref{eq:a1}. Then, knowing $b$ and $c$ allows us to determine $d:=\frac{N_1^2+N_2^2+N_3^2}{N_1 N_2 N_3}=b-2c$ and so we again have the integers  $p$, $q$, $r$ and $s$ from the proof of Theorem~\ref{thm:1}  which allow us to determine $N_1$, $N_2$ and $N_3$.  In summary, we have the following theorem.

\begin{Theorem}
\label{thm:3}
Let $M:=\C P^2(N_1,N_2,N_3)$ be a four-dimensional weighted projective space with isolated singularities.  If we fix the Euler characteristic of $M$ (i.e., fix the topological class), then the spectrum of the Laplacian acting on functions on $M$ determines the weights $N_1, N_2, N_3$.  
\end{Theorem}

\subsection{Hearing extremal metrics}
After determining the values of $N_1$, $N_2$ and $N_3$ we can also hear if a certain metric is extremal (cf. \S \ref{sec:A3}) from the value of $\mathfrak{a}_2$.

\begin{Theorem}
\label{thm:4}
Let $M:=\C P^2(N_1,N_2,N_3)$ be a four-dimensional weighted projective space with isolated singularities.  Suppose that the values of $N_1$, $N_2$, and $N_3$ are known.  Then the value of $\mathfrak{a}_2$ for the Laplacian acting on functions determines whether $M$ is endowed with the extremal metric.
\end{Theorem}

\begin{proof}
As shown in the proof of Theorem \ref{thm:1} (cf. \eqref{eqn:matrix}), 
\begin{equation}\label{eqn:a_2(0)}
\mathfrak{a}_2(\Delta_0) = \frac{1}{960 \pi^2} \int_M \scal^2 +  \frac{1}{90} \int_M  c_2 - \frac{1}{120} \int_M  c_1^2 \,\, + T.
\end{equation}
We have explicit expressions for $\int_M c_2$, $\int_M c_1^2$, and $T$ in terms of $N_1$, $N_2$, and $N_3$ (see \eqref{eq:b}, \eqref{eq:c}), 
so using \eqref{eqn:a_2(0)} we see that we can hear $\int_M \scal^2$ from the value of $\mathfrak{a}_2$.  If the metric is extremal, we have seen in \S \ref{sec:A3} that
\begin{equation}\label{eqn:scal2}
\int_M \scal^2 = 96 \pi^2 \frac{N_1^2 + N_2^2 + N_3^2}{N_1 N_2 N_3}.
\end{equation}
Hence, by checking the value we hear for $\int_M \scal^2$ against that given by \eqref{eqn:scal2}, we can see whether $M$ is endowed with the extremal metric.  Note that we are using the fact that the extremal metric is an absolute minimum for the $L^2$-norm of the scalar curvature among K\"ahler metrics in a fixed K\"ahler class.
\end{proof}



\begin{thebibliography}{99}

\bibitem{Ab} M. Abreu, {\em K\"ahler metrics on toric orbifolds},
  J. Differential Geom. {\bf 58} (2001), 151--187. 

\bibitem{A-C-G} V. Apostolov, D. Calderbank and P. Gauduchon, {\em Hamiltonian $2$-forms in K\"ahler geometry, I}, J. Differential Geom. {\bf 73} (2006), 359-412. 

\bibitem{At} M. Atiyah, {\em Convexity and commuting Hamiltonians}, Bull. London Math. Soc. {\bf 14} (1982), 1-15.


\bibitem{AB} M. Atiyah and R. Bott, {\em The moment map and equivariant cohomology}, Topology {\bf 23} (1984), 1-28.

\bibitem{BGM} Marcel Berger, Paul Gauduchon, and Edmond Mazet, \emph{Le
spectre d'une vari\'{e}t\'{e} riemannienne}, Lecture Notes in Mathematics, vol.
194, Springer-Verlag (1971).

\bibitem{BV} N. Berline and M. Vergne, {\em Classes caract\'eristiques \'equivariantes. Formule de localisation en cohomologie \'equivariante}, C. R. Acad. Sci. Paris S\'er. I Math. {\bf 295} (1982), 539-541.

\bibitem{B-Y} B. Berndt and B. Yeap, {\em Explicit evaluations and reciprocity theorems for finite trigonometric sums},  Adv. in Appl. Math. {\bf 29}  (2002), 358--385.  	

\bibitem{B} A. Besse, {\em Einstein manifolds}, Ergebnisse der Mathematik und ihrer Grenzgebiete, {\bf 10} Springer-Verlag, Berlin-New York, 1987.         

\bibitem{Br} R. Bryant, {\em Bochner-K\"ahler metrics},
J. Amer. Math. Soc. {\bf 14} (2001), 623--715.

\bibitem{Ca} E. Calabi, {\em Extremal K\"{a}hler metrics}, in 
``Seminar on Differential Geometry'' (ed. S.T.Yau), Annals of Math. Studies 
{\bf 102}, 259--290, Princeton Univ. Press, 1982.

\bibitem{Don76} H. Donnelly,  {\em Spectrum and the Fixed Point Sets of
Isometries I}, Math. Ann. {\bf 224} (1976), 161-170.

\bibitem{Don79} H. Donnelly, {\em Asymptotic expansions for the compact quotients of properly discontinuous group actions}, Illinois J. of Math. {\bf 23} (1979), 485-496.

\bibitem{DonPa} H. Donnelly and V.K. Patodi, {\em Spectrum and the fixed point sets of isometries. {II}}, Topology, {\bf 16} (1977), 1-11.

\bibitem{DGGWZ} E. Dryden, C. Gordon, S. Greenwald, D. Webb and S. Zhu, {\em Asymptotic Expansion of the Heat Kernel for Orbifolds}, preprint, 2006.

\bibitem{Far01} Carla Farsi, \emph{Orbifold Spectral Theory}, Rocky
Mountain J. Math., {\bf 31} (2001), 215-235.

\bibitem{FS} M. Furuta and B. Steer, {\em Seifert fibered homology $3$-spheres and the Yang-Mills equations on Riemann surfaces with marked points}, Adv. Math. {\bf 96} (1992), 38-102.

\bibitem{Gau} P. Gauduchon, {\em The Bochner tensor of a weakly Bochner-flat K\"{a}hler manifold}, An. Univ. Timi\c soara Ser. Mat.-Inform. {\bf 39} (2001), 189-237.

\bibitem{Gil04} P. Gilkey, {\em Asymptotic formulae in spectral geometry}, Studies in Advanced Mathematics, CRC Press, Boca Raton, FL, 2004. 

\bibitem{Go} L. Godinho, {\em Classification of symplectic circle actions on $4$-orbifolds}, in preparation.


\bibitem{H} A. Hwang, {\em On the Calabi energy of extremal K\"ahler metrics}, Internat. J. Math. {\bf 6} (1995), 825-830.


\bibitem{M} E. Meinrenken, {\em Symplectic surgery and the $Spin^c$-Dirac operator}, Adv. Math. {\bf 134} (1998), 240-277.

\bibitem{MOY} T. Mrowka, P. Ozsv\'{a}th and B. Yu {\em Seiberg-Witten Monopoles on Seifert Fibered Spaces}, Comm. Anal. Geom. {\bf 5} (1997), 685-791. 

\bibitem{R} S. Roan, {\em Picard groups of hypersurfaces in toric varieties}, Publ. Res. Inst. Math. Sci. {\bf 32} (1996), 797-834.

\bibitem{Sakai} T. Sakai, {\em On eigen-values of Laplacian and curvature of Riemannian manifold}, T\^ohoku Math. J. {\bf 23} (1971), 589-603.

\bibitem{S} P. Scott, {\em The geometries of $3$-manifolds}, Bull. London Math. Soc. {\bf 15} (1983), 401-487.

\bibitem{Se} H. Seifert, {\em Topologie dreidimensionaler gefaserter R\"aume}, Acta. Math. {\bf 60} (1932), 147-283. 
 
\bibitem{SSW} Naveed Shams, Elizabeth Stanhope, and David L. Webb,
\emph{One Cannot Hear Orbifold Isotropy Type}, preprint, 2005.

\bibitem{Sta05} Elizabeth A. Stanhope, \emph{Spectral Bounds on Orbifold
Isotropy}, Annals of Global Analysis and Geometry {\bf 27} (2005), 355-375. 
 
\bibitem{Ta1} S. Tanno, {\em Eigenvalues of the Laplacian of Riemannian manifolds}, T\^ohoku Math. J. {\bf 25} (1973), 391-403.

\bibitem{Ta2} S. Tanno, {\em An inequality for $4$-dimensional K\"ahlerian manifolds}, Proc. Japan Acad. {\bf 49} (1973), 257-261.
 
\bibitem{Z} D. Zagier, {\em Higher dimensional Dedekind sums},  Math. Ann.  {\bf 202}  (1973), 149--172. 


\end{thebibliography}
\end{document}